\documentstyle[11pt]{article}

\begin{document}

\noindent

\begin{center}

  {\LARGE Cohomology Ring of Crepant Resolutions of Orbifolds}

  \end{center}

  \noindent

  \begin{center}

    {\large  Yongbin Ruan}\footnote{ partially supported by the National Science Foundation and Hong Kong RGC}\\

    Department of Mathematics, Hong Kong University of Science and Technology\\

       and Department of Mathematics, University of Wisconsin-Madison\\

              \end{center}

              \def \x{{\bf x}}

              \def \M{{\cal M}}

              \def \C{{\bf C}}

              \def \Z{{\bf Z}}

              \def \R{{\bf R}}

              \def \Q{{\bf Q}}

              \def \U{{\cal U}}

              \def \E{{\cal E}}

              \def \z{{\bf z}}

              \def \m{{\bf m}}

              \def \n{{\bf n}}

              \def \g{{\bf g}}

              \def \h{{\bf h}}

              \def \V{{\cal V}}

              \def \W{{\cal W}}

              \def \T{{\cal T}}

              \def \P{{\bf P}}

        \def \L{{\cal L}}
    \def \H{{\cal H}}

    \section{Introduction}
    Suppose that $X$ is an orbifold. In general, $K_X$ is an orbifold vector bundle or a $Q$-divisor only.
    When the $X$ is so called  Gorenstein, $K_X$ is a bundle or a divisor. For Gorenstein orbifold, a
    resolution $\pi:Y \rightarrow X$ is called a crepant resolution if $\pi^*K_X=K_Y$. Here, "crepant"
    can be viewed as a minimality condition with respect to canonical bundle. Crepant resolution always
    exists when dimension is two or three. A nice way to construct it is to use Hilbert scheme of points.
    However, the crepant resolution in dimension three is not unique. Different crepant resolutions are
    connected by flops. When the dimension is bigger than four, the crepant resolution does not always exist.
    It is an extremely interesting problem in algebraic geometry to find out when it does exist. One of
    famous example is Hilbert scheme of points of algebraic surfaces, which is a crepant resolution of the symmetry
    product of an algebraic surface. When orbifold string theory was first constructed over global
quotient \cite{DHVW}, one of the first invariants from it is
orbifold Euler characteristic. It was conjectured that the
orbifold Euler characteristic is the same as the Euler
characteristic of its crepant resolution. This fits well with
McKay correspondence in algebro-geometry. It had been the main
attraction before the current development. By the work of Batyrev
and others \cite{B1} \cite{DL}, this conjecture has been extended
and solved for orbifold Hodge number of Gorenstein global
quotients. Very recently, it was solved in the complete
generalities by Lupercio-Poddar \cite{LP} and Yusuda \cite{Y}.

It has long been an interesting problem
    to compute cohomology of Hilbert scheme of points of algebraic surfaces. A great deal of works have been done
on this topic.  Among all the progress, a less known but crucial ingredient is the link to the orbifold cohomology
 of the symmetry     product.
  Suppose that $M$ is an algebraic surface. We use $M^{[n]}$ to denote the Hilbert scheme of points of length
    $n$ of $M$. In his thesis \cite{G}, G\"{o}ttsche computed the generating function of Euler number $\sum_{n=1}^{\infty} \chi(M^{[n]})q^n$
and showed that it has a surprising modularity. In 1994, in order
to explain its modularity, Vafa-Witten \cite{VW} computed
$\H=\oplus_n H^*(M^n/S_n,\C)$. Motivated by
    orbifold conformal field theory, they directly wrote $\H$ as a "Fock space" or a representation
    of Heisenberg algebra. Then, the generating function of Euler characteristic is interpreted as the correlation function of
    an elliptic curve. Therefore, it should be invariant under modular transformations of the elliptic curve. This shows that
the space of cohomology itself has more structure. Orbifold string
    theory conjecture predicates that $\oplus_n H^*(M^{[n]}, \C)$ should also admit a representation of Heisenberg algebra.
    This conjecture was verified by a beautiful work of Nakajima \cite{N} and others. One of theme of this short note is
that the orbifold  $M^n/S_n$ will continue to play a crucial role to compute ring structure of $H^*(M^{[n]}, \C)$.

    During last two years, there was a surge of activities to study mathematics of the orbifold string theory, which author called
    the stringy  orbifolds. A few curious physical concepts such as the orbifold Euler-Hodge
    numbers of global quotients found their places in a much broader and deeper theory. For example, a new cohomology (orbifold cohomology)
    was constructed \cite{CR1}. The growth was so explosive that the author believes
    that there is an emerging new subject of mathematics. He learned in graduate school that the test of the relevance of a new theory
    has been the progress it made from old problems. Therefore, it is particularly significant to revisit the problem of computing
    the ring structure of $M^{[n]}$ and a crepant resolution in general. A lot of information is known for additive
structure of $H^*(M^{[n]}, \C)$. The ring structure of $M^{[n]}$ is quite subtle and more
    interesting. Partial results has been obtained by
    Fantechi-G\"{o}ttsche \cite{FG1}, Ellingsrud-Stromme
    \cite{ES1} \cite{ES2}, Beauvill \cite{Bea}, Mark \cite{Mar}.
Based on  an important observation by Frenkel
    and Wang \cite{FW}, Lehn-Sorger \cite{LS1} determined the
    cohomology ring of $(\C^2)^{[n]}$. At the same time, the
    author was computing orbifold cohomology $(\C^2)^n/S_n$ and
    the result  from both calculations matches perfectly.  Based on physical motivation and the strong evidence from $(\C^2)^n/S_n$.
    The author proposed \cite{R2} a conjecture in the case
    of a hyperkahler resolution.

    \vskip 0.1in
    \noindent
    {\bf Cohomological Hyperkahler Resolution Conjecture: }{\it Suppose that $\pi: Y\rightarrow X$ is a hyperkahler resolution. Then, the ordinary cohomology ring of $Y$
    is isomorphic to the orbifold cohomology ring of $X$.}
    \vskip 0.1in
    In the case of Hilbert scheme points of surfaces, the Cohomological Hyperkahler Resolution Conjecture (CHRC) implies that $K3^{[n]}, (T^4)^{[n]}$ have isomorphic cohomology ring
    as the orbifold cohomology rings of $K3^n/S_n, (T^4)^{n}/S_n$. The later was proved recently by the beautiful works
    of Lehn-Sorge \cite{LS2}, Fantechi-G\"{o}ttsche \cite{FG2} and Uribe \cite{U}.
    It should be mentioned that  Fantechi-G\"{o}ttsche-Uribe's work computed the orbifold ring structure of $X^n/S_n$ for an arbitrary
    complex manifold $X$ which may or may not be $K3, T^4$. There is a curious phenomenon that over rational number Lehn-Sorge, Fantanch-G\"{o}ttsche and Uribe showed that one must modify the ring structure of the orbifold
    cohomology by a sign in order to match to the cohomology of Hilbert scheme. However, Qin-Wang observed that
	such a sign modification is unnecessary over complex number \cite{QW}. All the conjectures stated in this article (in fact
	any conjecture motivated by physics) are the statements over complex number. 
     The ring structure of $X^{[n]}$ for a general algebraic surface $X$ is
    still unknown.

    It is easy to check that CHRC is false if we drop the hyperkahler condition. One of main purposes of this article is to
    propose a conjecture for the arbitrary crepant resolution.

    As  mentioned previously, the crepant resolutions are not unique. The different crepant resolutions are connected by
    "K-equivalence" \cite{W}.  Two smooth (or Gorenstein orbifolds) complex
    manifolds $X, Y$ are $K$-equivalent iff there is a common resolution $\phi, \psi: Z\rightarrow X, Y$ such that $\phi^* K_X=
    \psi^* K_Y$. Batyrev-Wang \cite{B}, \cite{W} showed that two $K$-equivalent projective manifolds have the same betti number. It is natural to
    ask if they have the same ring structures. This question is obviously related to CHRC. Suppose that CHRC holds for non-hyperhahler
    resolutions. It implies that different resolutions (K-equivalent) have the same ring structures. Unfortunately, they usually have
    different ring structures,  and hence CHRC fails in general. It is
easy to check this
    in case of  three dimensional flops. A key idea to remedy the  situation is to include the quantum corrections. The author proposed \cite{R1}
    \vskip 0.1in
    {\bf Quantum Minimal Model Conjecture: }{\it Two $K$-equivalent projective manifolds have the same quantum cohomologies.}
    \vskip 0.1in
    Li and the author proved Quantum Minimal Model conjecture in complex dimension three. In higher dimensions, it seems to be
    a difficult problem. In many ways, Quantum Minimal Model Conjecture unveils the deep relation between the quantum cohomology and
    the birational geometry \cite{R2}. However, it is a formidable task to master the quantum cohomology machinery for any non-experts.
     In this article, we proposed another conjecture focusing on the cohomology
    instead of the quantum cohomology. As mentioned before, the cohomology ring structures are not isomorphic for $K$-equivalent
    manifolds. Therefore, some quantum information must be included. Our new conjecture requires a minimal set of quantum information
    involving the GW-invariants of exceptional
    rational curve.

    Finally,  the motivation behind the conjectures should be described here. Let's first go back to the motivation of CHRC. Naively,
    physics indicates that the orbifold quantum cohomology of $X$ should be "equivalent" to the quantum cohomology of $Y$. It is not clear how to
      formulate the precise meaning of the "equivalence". However, for the hyperkahler resolutions, there are no quantum corrections and
     the quantum     cohomology is just cohomology. All the difficulties to formulate the "equivalence" of the quantum information disappear. We should just
    get an isomorphism between cohomologies. This is the reasoning behind Cohomological Hyperkahler Resolution Conjecture. For a general
    crepant resolution, quantum corrections do appear. Mathematically, it means that the cup product of a crepant resolution is the orbifold
    cup product of the orbifold plus some quantum corrections. A further study shows that the quantum corrections  come from the GW-invariants
    of exceptional rational curves only. However, when we try to count these quantum corrections, we encounter a serious problem. The
    quantum corrections appear to be an infinite series of the GW-invariants corresponding to the multiple degree of exceptional rational
    curves. In the quantum cohomology, we insert a quantum variable $q$ to keep track of this degree. Intuitively, we should set $q=1$. By checking
    a few examples, one can find that the quantum corrections diverge at $q=1$ and some kind of "renormalization" is necessary.
    If we believe that a solution can be found in the other value of
    $q$, it is instructive to find the possible value of $q$. When
    we match the quantum cohomologies under 3-dimensional flops, there
    is a change of quantum variable $q\rightarrow \frac{1}{q}$. If
    there is an uniform way to set the value of $q$, said
    $q=\lambda$. Then, we must have $\lambda=\frac{1}{\lambda}$.
    Hence, $\lambda^2=1$. Namely, the other choice is $q=-1$.
    In a different context, there was  beautiful works by P. Aspinwall \cite{A} and E. Wendland \cite{W} in physics  concerning the conformal theory on $K3$ \cite{W}.
    It suggests that after the quantization, the value of B-field "shifts" to
    "$q=-1$"!

    This short article is organized as follows. In the section two, we will formulate our conjectures. In the section three, we will
    verify our conjectures using several examples. Special thanks goes to P. Aspinwall and E. Witten for bringing  me the attention of
    \cite{W} and K. Wendland for a wonderful talk on Workshop on Mathematical Aspect of Orbifold String Theory to explain \cite{W}.
    Finally, I would like to thank Wei-Ping Li, Zhenbo Qin  for interesting discussions.

    \section{Conjectures}
       Suppose that $\pi: Y\rightarrow X$ is one crepant resolution of Gorenstein  orbifold $X$.  Then, $\pi$ is
	a Mori contraction and the homology classes of rational curves $\pi$ contracted are generated by so called
	extremal rays. Let $A_1, \cdots,
    A_k$ be an integral basis of extremal rays. We call $\pi$ non-degenerate if  $A_1, \cdots,
    A_k$ are linearly independent. For example, the Hilbert-Chow
	map $\pi: M^{[n]}\rightarrow M^n/S_n$ satisfies this hypothesis.  Then,
    the homology class of any effective curve being contracted can be written as
    $A=\sum_i a_i A_i$ for $a_i\geq 0$. For each $A_i$, we assign a formal variable $q_i$. Then, $A$ corresponds to $q^{a_1}_1 \cdots q^{a_k}_k$.
    We define a 3-point function
    $$<\alpha, \beta, \gamma>_{qc}(q_1, \cdots q_k)=\sum_{a_1, \cdots, a_k}\Psi^X_{A}(\alpha, \beta, \gamma)q^{a_1}_1\cdots
    q^{a_k}_k ,\leqno(2.1)$$
    where $\Psi^X_{A}(\alpha, \beta, \gamma)$ is Gromov-Witten invariant and $qc$ stands for the quantum correction.
    We view $<\alpha, \beta, \gamma>_{qc}(q_1, \cdots, q_k)$
    as  analytic function of $q_1, \cdots q_k$ and  set $q_i=-1$ and let
    $$<\alpha, \beta,\gamma>_{qc}=<\alpha, \beta,\gamma>_{qc}(-1, \cdots, -1).\leqno(2.2)$$
    We define a quantum corrected triple intersection
    $$<\alpha, \beta, \gamma>_{\pi}=<\alpha, \beta,
    \gamma>+<\alpha, \beta, \gamma>_{qc},$$
    where $<\alpha, \beta, \gamma>=\int_X \alpha\cup\beta\cup\gamma$ is the ordinary triple intersection.
    Then we define the quantum corrected cup product $\alpha\cup_{\pi}
    \beta$ by the equation
    $$<\alpha\cup_{\pi} \beta, \gamma>=<\alpha, \beta,
    \gamma>_{\pi},$$
    for arbitrary $\gamma$. Another way to understand $\alpha\cup_{\pi}\beta$ is as following.
    Define a product $\alpha\star_{qc}\beta$ by the equation
    $$<\alpha\star_{qc}\beta, \gamma>=<\alpha, \beta,\gamma>_{qc}$$
    for arbitrary $\gamma \in H^*(Y, \C)$. Then, the quantum corrected product is the ordinary cup product
    corrected by $\alpha\star_{qc}\beta$. Namely,
    $$\alpha\cup_{\pi} \beta=\alpha\cup \beta +\alpha\star_{qc} \beta. \leqno(2.3) $$
    We denote the new quantum corrected cohomology ring as $H^*_{\pi}(Y, \C)$.

    \vskip 0.1in
    \noindent
    {\bf Cohomological Crepant Resolution Conjecture: }{\it Suppose that $\pi$ is non-degenerate and hence
	$H^*_{\pi}(Y, \C)$ is well-defined. Then,  $H^*_{\pi}(Y, \C)$
              is the ring isomorphic to orbifold cohomology ring $H^*_{orb}(X, \C)$.}
    \vskip 0.1in
    Recently, Li-Qin-Wang \cite{LQW3} proved a striking theorem that the
    cohomology ring of $M^{[n]}$ is universal in the sense that it
    depends only on homotopy type of $M$ and $K_X$. Combined with
    their result,  CCRC yields
    \vskip 0.1in
    \noindent
    {\bf Conjecture: }{\it For $M^{[n]}$,  $\star$ product depends only on $K_M$.}
    \vskip 0.1in
    It suggests an interesting way to calculate $\star$ product by
    first finding a universal formula (depending only on $K_M$)
    and calculating a special example such as $\P^2$ to determine
    the coefficient.

    Next, we formulate a closely related conjecture for $K$-equivalent manifolds.

    Suppose that $X, X'$ are $K$-equivalent and $\pi: X\rightarrow X'$ is the birational map. Again, exceptional rational curves makes sense.
     Suppose that $\pi$ is nondegenerate. Then, we go through the previous
    construction to define ring $H^*_{\pi}(X, \C)$.
    \vskip 0.1in
    \noindent
    {\bf Cohomological Minimal Model Conjecture: }{\it Suppose that $\pi, \pi^{-1}$ are nondegenerate. Then,
    $H^*_{\pi}(X, \C)$ is the ring isomorphic to $H^*_{\pi^{-1}}(X', \C)$ }
    \vskip 0.1in
    Here $\pi^{-1}: X'\rightarrow X$ is the inverse birational transformation of $\pi$.
    When $X, X'$ are the different crepant resolutions of the same orbifolds, Cohomological minimal model conjecture follows from Cohomological
    crepant resolution conjecture. However, it is well-known that most of K-equivalent manifolds are not crepant resolution of orbifolds.
    Cohomological minimal model conjecture can be generalized to orbifold provided that the quantum corrections are defined using orbifold
    Gromov-Witten invariants introduced by Chen-Ruan \cite{CR2}. 
\vskip 0.1in
    \noindent
    {\bf Remark 3.4: }{\it (1) The author does not know how to define quantum corrected cohomology if $\pi$ is
not nondegenerate. (2) All the conjectures in this section
    should be understood as the conjectures up to certain slight
    modifications (see next sections).}
\vskip 0.1in

    \section{Verification of Conjectures}
        \vskip 0.1in
    \noindent
    {\bf Example 3.1: } Suppose that $\Sigma$ is one Riemann surface of genus $\geq 2$ and $E\rightarrow \Sigma$ is a rank two
    bundle  such that $C_1(E)=2g-2$. Then, $E$ is an example of local Calabi-Yau manifold. Let $\tau$ be the involution acting
    on $E$ as the multiplication of $-1$. $X=E/\tau$ is a Calabi-Yau orbifold. Let $\tilde{E}$ be the blow-up of $E$ along $\Sigma$. The
    action of $\tau$ extends over $\tilde{E}$. Let $Y=\tilde{E}/\tau$. The projection $\pi:Y \rightarrow X$ is a crepant resolution of
    $X$. Let's verify Cohomological Crepant Resolution Conjecture in this case. For simplicity, we consider the even cohomology
    only. Moreover, it is enough to compare triple intersections $<\alpha, \beta, \gamma>$ for $\alpha, \beta, \gamma\in H^2$.

    Note that $X$ is homotopic equivalent to $\Sigma$. Therefore, the nontwisted sector contributes one generator to $H^2_{orb}(X, \C)$.
    Let $\alpha$ be the generator with the integral one on $\Sigma$. Since $\Sigma$ has local group $\Z_2$,
      it generates a twisted sector with    degree shifting number $1$. It contributes to a generator $\beta$ to $H^2_{orb}(X, \C)$, where
    $\beta$ represents the  constant function $1$ on the twisted sector. It is easy to compute
    $$<\alpha, \alpha, \alpha>=0, <\alpha, \alpha, \beta>=0, <\alpha, \beta, \beta>=\frac{1}{2}, <\beta, \beta, \beta>=0.\leqno(3.1)$$

    Let's compute $H^{2}(Y, \C)$. Let $\alpha'=\pi^*\alpha$. The exceptional divisor $S$ is a ruled surface of $\Sigma$.
      Let $\beta'$ be its Poincare dual. It is clear that
    $$<\alpha',\alpha', \alpha'>=0, <\alpha', \alpha',\beta'>=0.\leqno(3.2)$$
    Nonzero ones are $<\alpha',\beta',\beta'>=K_S[C]=-2, <\beta',\beta',\beta'>=K^2_S=8(1-g)$ \cite{W} (Lemma 3.2),
    where $K_S$ is the canonical bundle of $S$ and $C$ is the fiber of $S$.

    Next, we compute the correction term. By Wilson \cite{W}(Lemma 3.3), a small complex deformation can deform $S$ into $2(g-1)$
    many $\O(-1)+\O(-1)$ curves. Wilson's argument is completely local and works for this case. It is well-known that
    $$<\gamma_1, \gamma_2, \gamma_3>_{qc}(q)=2(g-1)\gamma_1([C])\gamma_2([C])\gamma_3([C])\sum_{i=1} q^i=
    2(g-1)\gamma_1([C])\gamma_2([C])\gamma_3([C])\frac{q}{1-q}.\leqno(3.3)$$
    Hence,
    $$<\gamma_1, \gamma_2, \gamma_3>_{qc}=(1-g)\gamma_1([C])\gamma_2([C])\gamma_3([C]).\leqno(3.4)$$
    Moreover, $\alpha'[C]=0, \beta'[C]=K_S[C]=-2$. Then, we obtain that
    $$<\alpha',\alpha',\alpha'>_{qc}=<\alpha', \alpha', \beta'>_{qc}=<\alpha',\beta', \beta'>_{qc}=0.\leqno(3.5)$$
    $$<\beta',\beta',\beta'>_{qc}=-8(1-g).\leqno(3.6)$$
    After corrected by $<,,,>_{qc}$. $<\beta, \beta, \beta>, <\beta', \beta', \beta'>_{\pi}$ match perfectly. But there is still a discrepancy
    between $<\alpha, \beta, \beta>, <\alpha', \beta', \beta'>_{\pi}$ which is the reminisce of discrepancy for surface quotient singularities
    (See Remark). We don't  know a canonical way to construct a homomorphism between orbifold cohomology and cohomology of its crepant resolution.
    Nevertheless, up to a sign, the map $\alpha\rightarrow \alpha', \beta\rightarrow 2\beta'$ gives an ring isomorphism. $\Box$

    \vskip 0.1in
    \noindent
    {\bf Example 3.2: } Next, we use the work of Li-Qin \cite{LQ}
    to verify Cohomological Crepant Resolution Conjecture for
    $M^{[2]}$. To simplify the formula, we assume that $M$ is simply
    connected.

    It is easy to compute the orbifold cohomology $H^*_{orb}(X,
    \C)$ for $X=M^2/\Z_2$. The nontwisted sector
    can be identified with invariant cohomology of $M^2$.
    Let $h_i\in H^2(M, \C)$ be a basis and $H\in H^4(M, \C)$ be Poincare dual to
    a point. Then, the cohomology of the nontwisted
    sectors are generated by $1, 1\otimes h_i+h_i\otimes 1, 1\otimes
    H+H\otimes 1, h_i\otimes h_j+h_j\otimes h_i, h_i\otimes H+H\otimes h_i,
    H\otimes H.$ The twisted sector is diffeomorphic to
    $M$ with degree shifting number 1. We use $\tilde{1},
    \tilde{h_i}, \tilde{H}$ to denote the generators. They are of
    degrees $2,4,6$. By the definition, triple intersections
    $$<twisted sector, nontwisted sector, nontwisted sector>=0,$$
    $$<twisted sector, twisted sector, twisted sector>=0.$$
    Following is the table of nonzero triple intersections
    involving classes from the twisted sector
    $$<\tilde{1}, \tilde{1}, 1\otimes H+H\otimes 1>=1, <\tilde{1},
    \tilde{1}, h_i\otimes h_j+h_j\otimes h_i>=<h_i, h_j>, $$
    $$ <\tilde{1}, \tilde{h}_i, h_j\otimes 1+h_j\otimes 1>=<h_i, h_j>. \leqno(3.7)$$

    Next, we review the construction of $Y=M^{[2]}$.
    Let $\widetilde{M^2}$ be the blow-up of
    $M^2$ along the diagonal. Then, $\Z_2$ action
    extends to $\widetilde{M^2}$. Then, $Y=\widetilde{M^2}/\Z_2$.
    It is clear that we should map the classes from nontwisted sector to its pull-back $\pi: Y\rightarrow X$.
    We use the same notation to denote them. The exceptional divisor $E$ of Hilbert-Chow map
    $\pi: Y\rightarrow X$ is a $P^1$-bundle over $M$. Let
    $\bar{1}, \bar{h}_i, \bar{H}$ be the Poincare dual to $E$,
     $p^{-1}(PD(h_i))$ and fiber $[C]$, where $p: E\rightarrow M$ is the projection.

    Notes that $\bar{1}|_E=2\E$, where $\E$ is the tautological
    divisor of $P^1$-bundle $E\rightarrow M$. It is clear that
    $E=P(N_{\Delta(X)|X^2})$, where $\Delta(X)\subset X^2$ is the diagonal. Hence,
    $$<\bar{1}, \bar{1},
    \bar{h}_i>=4\E^2|_{p^{-1}(PD(h_i))}=4C_1(N_{\Delta(X)|X^2})\E|_{p^{-1}(PD(h_i))}=-4<C_1(X),
    h_i>. \leqno(3.8)$$
    $$<\bar{1}, \bar{1}, 1\otimes H+H\otimes 1>=2\E(1\otimes
    H+H\otimes 1)(E)=4\E(C)=-4.\leqno(3.9)$$
    $$<\bar{1}, \bar{1}, h_i\otimes h_j+h_j\otimes h_i>=-4<h_i,
    h_j>.\leqno(3.10)$$
    $$<\bar{1}, \bar{h}_i, 1\otimes h_j+h_j\otimes 1>=-4<h_i,
    h_j>.\leqno(3.11)$$
    Others are zero.

    The quantum corrections have been computed by Li-Qin
    \cite{LQ} (Proposition 3.021). The only nonzero terms are
    $$\begin{array}{lll}
    <\bar{1}, \bar{1}, \bar{h}_i>_{qc}(q)&=&\sum_{d=1}
    \bar{1}(d[C])^2\Psi^X_{d[C]}(\bar{h})q^d\\
    &=&\sum_{d=1} \frac{4d^2(2<K_X, h_i>)}{d^2} q^d\\
    &=&8<K_X, h_i>\frac{q}{1-q}
    \end{array}.\leqno(3.12)$$
    Hence,
    $$<\bar{1}, \bar{1}, \bar{h}_i>_{qc}=4<K_X, h_i>=4<C_1(X),
    h_i>\leqno(3.13)$$
    cancels $<\bar{1}, \bar{1},\bar{h}_i>$.

    It is clear that the map $\tilde{1}\rightarrow 2\bar{1},
    \tilde{h}\rightarrow 2\bar{h}, \tilde{H}\rightarrow \bar{H}$ is
    a ring isomorphism. $\Box$

    Next, we give two examples to verify  Cohomological Minimal Model Conjecture (CMMC).
    \vskip 0.1in
    \noindent
    {\bf Example 3.4: } The first example is the flop in dimension
    three. This case has been worked out in great detail by Li-Ruan
    \cite{LR}. For example, they proved a theorem that quantum cohomology rings are
    isomorphic under the change of the variable $q\rightarrow
    \frac{1}{q}$. Notes that if we set $q=-1$, $\frac{1}{q}=-1$. We set other quantum variables zero.
    Then, the quantum product becomes the quantum corrected product
    $\alpha\cup_{\pi}\beta$. Hence, CMMC follows from Li-Ruan's
    theorem. However, It should be pointed out that one can directly
    verify CMMC without using Li-Ruan's theorem. In fact, it is
    an much easier calculation.

    \vskip 0.1in
    \noindent
    {\bf Example 3.5: } There is a beautiful four dimensional
    birational
    transformation called Mukai transform as follows. Let
    $\P^2\subset X^4$ with $N_{\P^2|X^4}=T^* \P^2$. Then, one can
    blow up $\P^2$. The exceptional divisor of the blow up is  a hypersurface of $\P^2\times \P^2$
    with the bidegree $(1,1)$. Then, one can blow down in another
    direction to obtain $X'$. $X, X'$ are $K$-equivalent. In his
    Ph.D thesis \cite{Z}, Wanchuan Zhang showed that the quantum
    corrections $<\alpha, \beta, \gamma>_{qc}$ are trivial, and
    cohomologies of $X, X'$ are isomorphic.

    \section{Remarks}

    In the computation of orbifold cohomology of symmetry product
    and its relation to that of Hilbert scheme of points, there
    are two issues arisen. It was showed in the work of
    Lehn-Sorger, Fantechi-G\"{o}ttsche and Uribe that one has to
    add a sign in the definition of orbifold product in order to
    match that of Hilbert scheme of points over rational number. This sign was
    described as follows.

    Recall the definition of orbifold cup
    product
    $$\alpha\cup_{orb}\beta =\sum_{(h_1, h_2)\in T_2, h_i\in
    (g_i)} (\alpha\cup_{orb}\beta)_{(h_1,h_2)},\leqno(3.19)$$
    where $(\alpha\cup_{orb}\beta)_{(h_1,h_2)}\in
    H^*(X_{(h_1h_2)}, \C)$ is defined by the relation
    $$<(\alpha\cup_{orb}\beta)_{(h_1,h_2)}, \gamma>_{orb}
    =\int_{X_{(h_1,h_2)}}e^*_1\alpha\wedge e^*_2\beta\wedge e^*_3\gamma \wedge
    e_A(E_{(\g)}).\leqno(3.20)$$
    for $\gamma\in H^*_c(X_{((h_1h_2)^{-1})}, \C)$
    Then, we add a sign to each term.
    $$\alpha\cup_{orb}\beta =\sum_{(h_1, h_2)\in T_2, h_i\in
    (g_i)} (-1)^{\epsilon(h_1, h_2)}(\alpha\cup_{orb}\beta)_{(h_1,h_2)},\leqno(3.21)$$

    where
    $$\epsilon(h_1,
    h_2)=\frac{1}{2}(\iota(h_1)+\iota(h_2)-\iota(h_1h_2)).\leqno(3.22)$$
    Since
    $$\epsilon(h_1, h_2)+\epsilon(h_1h_2,
    h_3)=\frac{1}{2}(\iota(h_1)+\iota(h_2)+\iota(h_3)-\iota(h_1h_2h_3))=\epsilon(h_1,
    h_2h_3)+\epsilon(h_2, h_3),$$
    such a sign modification does not affect the associativity of
    orbifold cohomology.

	However, Qin-Wang \cite{QW} observed that the orbifold cohomology modified by such a sign
	is isomorphic to original orbifold cohomology over complex number by an explicit isomorphism
 	$$\alpha\rightarrow (-1)^{\frac{\iota(g)}{2}} \alpha$$
	for $\alpha\in H^*(X_{(g)}, \C)$. $\epsilon(h_1, h_2)$ is often an integer (for
	example symmetric product) while $\frac{\iota(g)}{2}$ is just a fraction. Hence, 
	$(-1)^{\frac{\iota(g)}{2}} $ is a complex number only.

    Another issue is the example of the crepant resolution of surface singularities
      $\C^2/\Gamma$. As Fantechi-G\"{o}ttsche \cite{FG2} pointed out, the Poincare paring of $H^2_{orb}(\C^2/\Gamma, \C)$
    is indefinite while the Poincare paring of its crepant resolution is negative definite. There is an
    easy way to fix this case (suggested to this author by Witten). We view the involution $I: H^*(X_{g^{-1}}, \C)\rightarrow H^*(X_g, \C)$
    as a "complex conjugation". Then, we define a "hermitian inner product"
    $$<<\alpha, \beta>>=<\alpha, I^*(\beta)>.\leqno(3.23)$$
    If we use this "hermitian" inner product, the intersection paring is positive definite again.
     The above process has its conformal theory origin (see
     \cite{NW}). It is attempting to perform this modification on
     $H^*(X_{(g)}, \C)\oplus H^*(X_{(g^{-1})}, \C)$ whenever
     $\iota_{(g)}=\iota_{(g^{-1})}$.
  The author does not know if it will affect  the associativity of
      orbifold cohomology.

\end{document}